\documentclass{amsart}

\usepackage{graphicx}
\usepackage[graphicx]{realboxes}
\newtheorem{thm}{Theorem}
\newtheorem{defn}{Definition}

\newtheorem{rem}{Remark}
\newtheorem{exa}{Example}
\newtheorem{conj}{Conjecture}

\usepackage[pdfauthor="Richard J. Mathar",colorlinks=true,bookmarks,citecolor=blue]{hyperref}

\begin{document}

\title[Non-bonding Dominoes]{Bivariate Generating Functions Enumerating Non-Bonding Dominoes on Rectangular Boards}

\author{Richard J. Mathar}
\urladdr{https://www.mpia-hd.mpg.de/~mathar}
\email{mathar@mpia-hd.mpg.de}
\address{Max-Planck Institute for Astronomy, K\"onigstuhl 17, 69117 Heidelberg, Germany}


\date{\today}

\begin{abstract}
The manuscript studies configurations of non-overlapping non-bonding
dominoes on finite rectangular boards of unit squares characterized 
by row and column
number. The non-bonding dominoes are defined here by the requirement
that any domino on the board shares at most one point (one
of its four corner points) with any other domino, but no edge.
With the Transfer Matrix Method, rational generating functions
are derived that solve the enumeration problem entirely, here evaluated
for boards with up to six  rows or columns.
\end{abstract}

\maketitle
\section{Non-Bonding Dominoes} 

A domino covers two adjacent squares on the square grid
vertically or horizontally. 
We call two dominoes non-bonding (or non-adjacent)
when they do not share any of their
7 (1 internal + 6 perimeter) edges, that is, if they share
at most one point at one of the four corners.
The same criterion is that the distance of any of the two 
squares in a domino has 
minimum L1 (Manhattan) distance of 2 to any other square
in a different domino.

In a raster-scan black-and-white images of a board partially filled with non-bonding dominoes,
the confusion limit is never passed: the dominoes remain distinct
and never can be mistaken as
tetrominoes, hexominoes etc.

\begin{rem}
The non-bonding monominoes have been studied by Siehler \cite{SiehlerArxiv1409,MatharVixra2404}.
\end{rem}

\begin{defn}
$D(r,c,d)$ is the number of arrangements of $d$ non-overlapping,
non-bonding dominoes
on a $r\times c$ rectangular square grid.
\end{defn}
The role of rows
and columns may be interchanged:
\begin{equation}
D(r,c,d)=D(c,r,d).
\label{eq.Dsym}
\end{equation}

\begin{figure}
\includegraphics[scale=0.3]{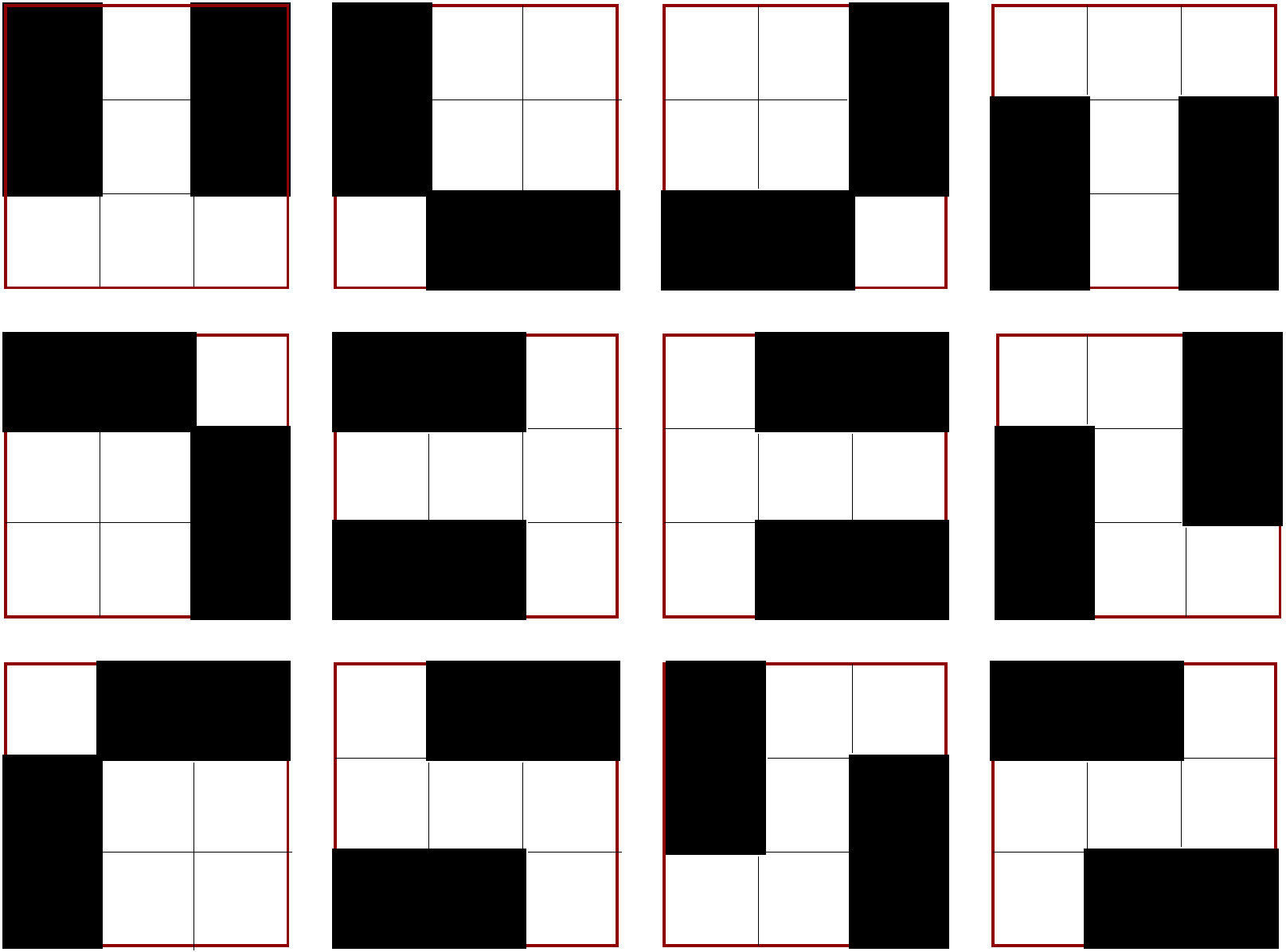}
\caption{All $D(3,3,2)=12$ configurations  of 2 non-bonding dominoes on the $3\times 3$ board.
}
\label{Fig.332}
\end{figure}

\begin{figure}
\includegraphics[scale=0.3]{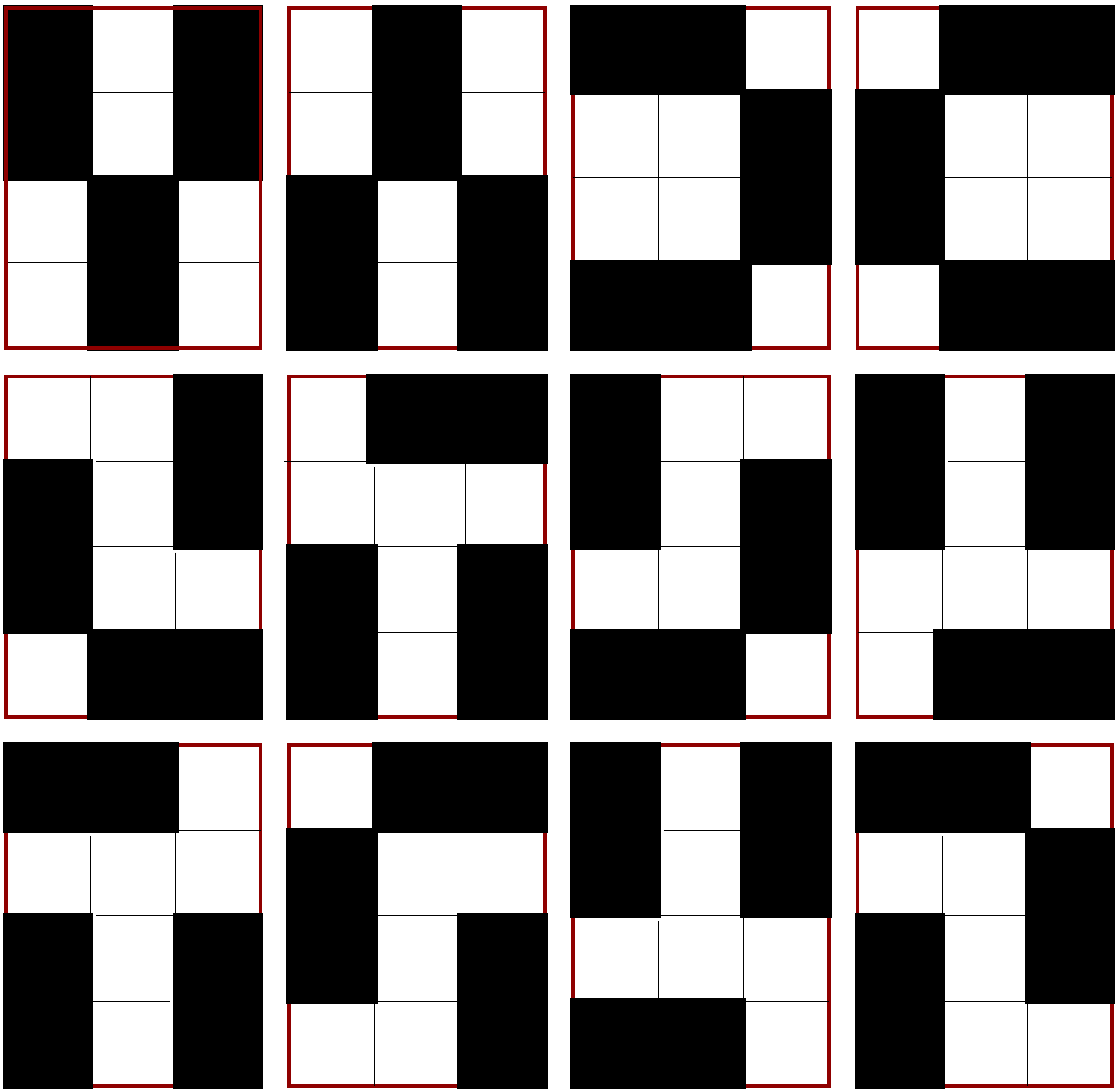}
\caption{All $D(4,3,3)=12$ configurations of 3 non-bonding dominoes on the $4\times 3$ board.
}
\label{Fig.433}
\end{figure}

Figures \ref{Fig.332} and \ref{Fig.433} illustrate the allowed arrangements
for the maximum filling of the $3\times 3$ and $4\times 3$ boards.
\begin{figure}
\includegraphics[scale=0.25]{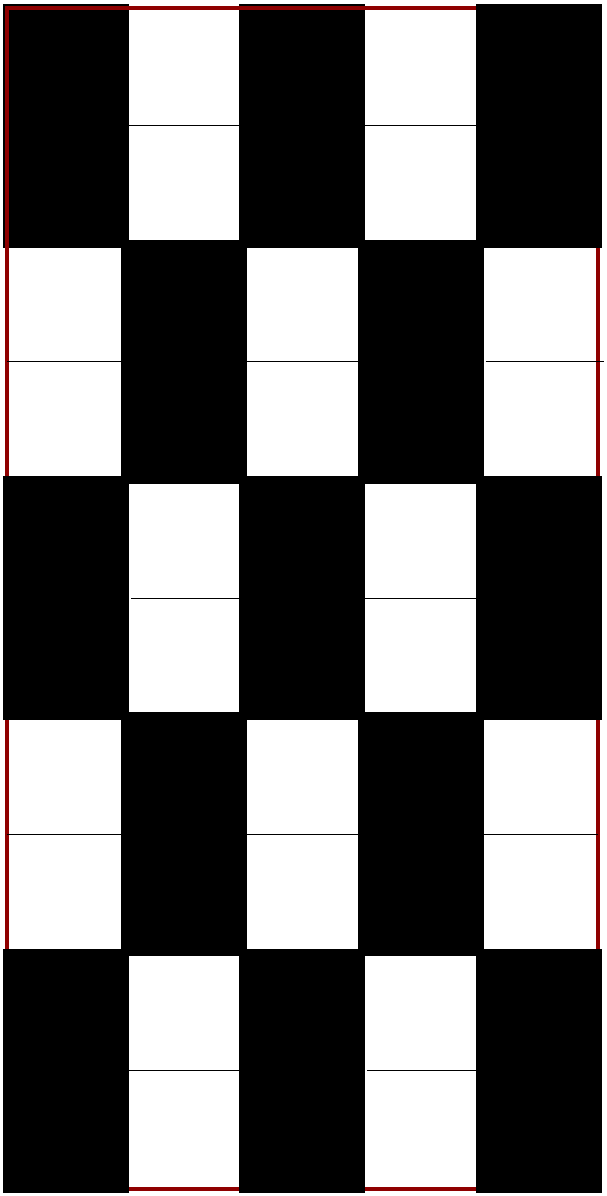}
\caption{The $D(10,5,13)=1$ way of placing 13 non-bonding dominoes on the $10\times 5$ board.
If $r=2(2r'+1)$ is an odd multiple of 2 and $c=2c'+1$ is odd (or vice versa), one can place 
$\bar d=(r'+1)(c'+1)+r'c'=rc/4+1/2$
dominoes adopting that regular pattern.
}
\label{Fig.10513}
\end{figure}

\begin{figure}
\includegraphics[scale=0.3]{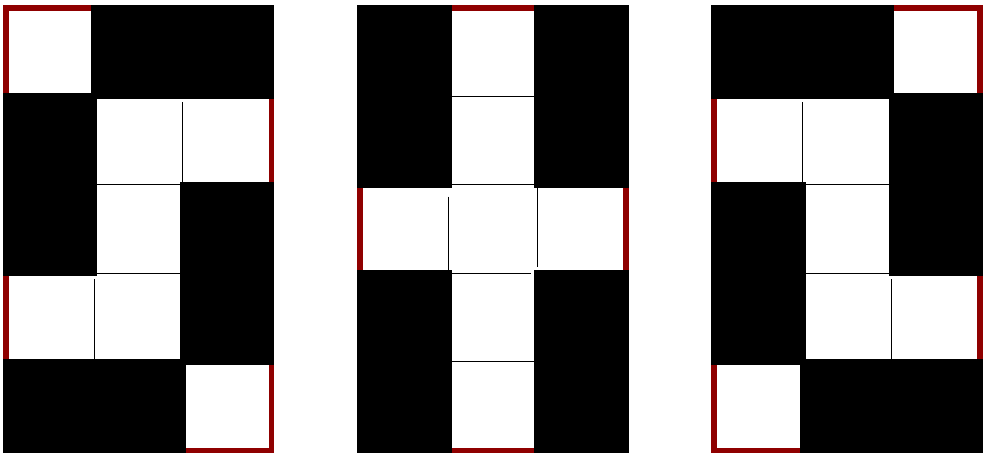}
\caption{The $D(5,3,4)=3$ configurations of 4 non-bonding dominoes on the $5\times 3$ board.
}
\label{Fig.534}
\end{figure}

\begin{figure}
\includegraphics[scale=0.4]{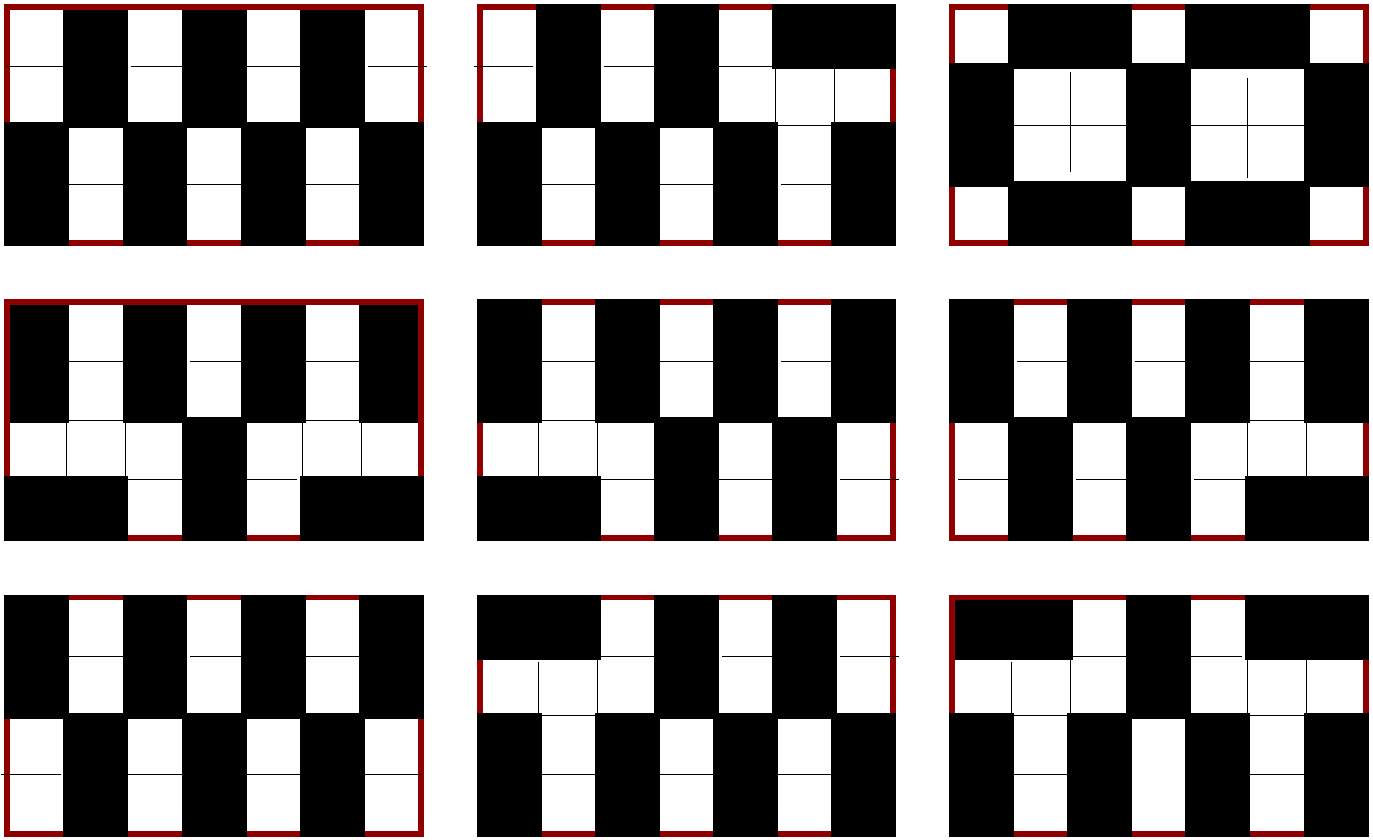}
\caption{The $D(7,4,7)=9$ configurations of 7 non-bonding dominoes on the $7\times 4$ board.
}
\label{Fig.747}
\end{figure}

\begin{figure}
\includegraphics[scale=0.4]{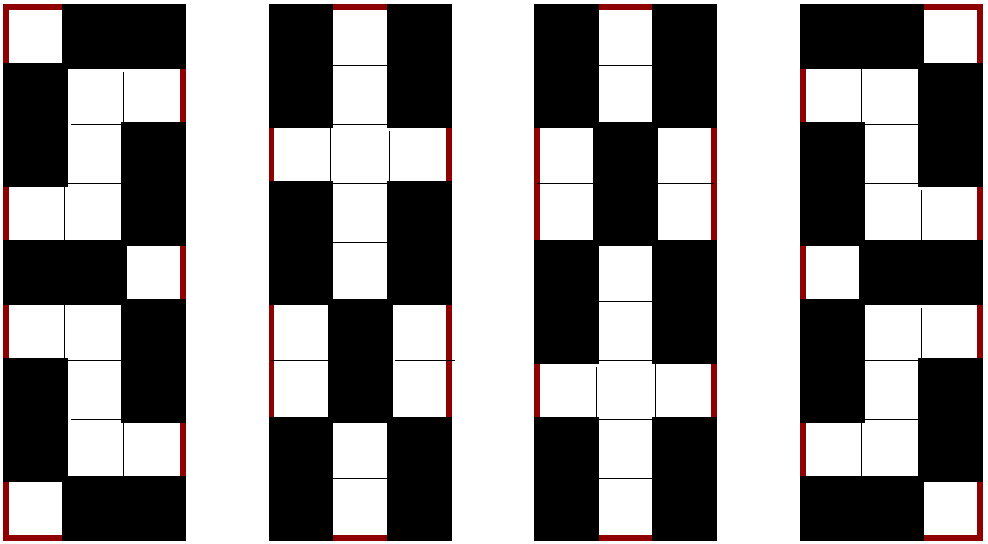}
\caption{The $D(9,3,7)=4$ configurations of 7 non-bonding dominoes on the $9\times 3$ board.
}
\label{Fig.937}
\end{figure}

\begin{figure}
\includegraphics[scale=0.5]{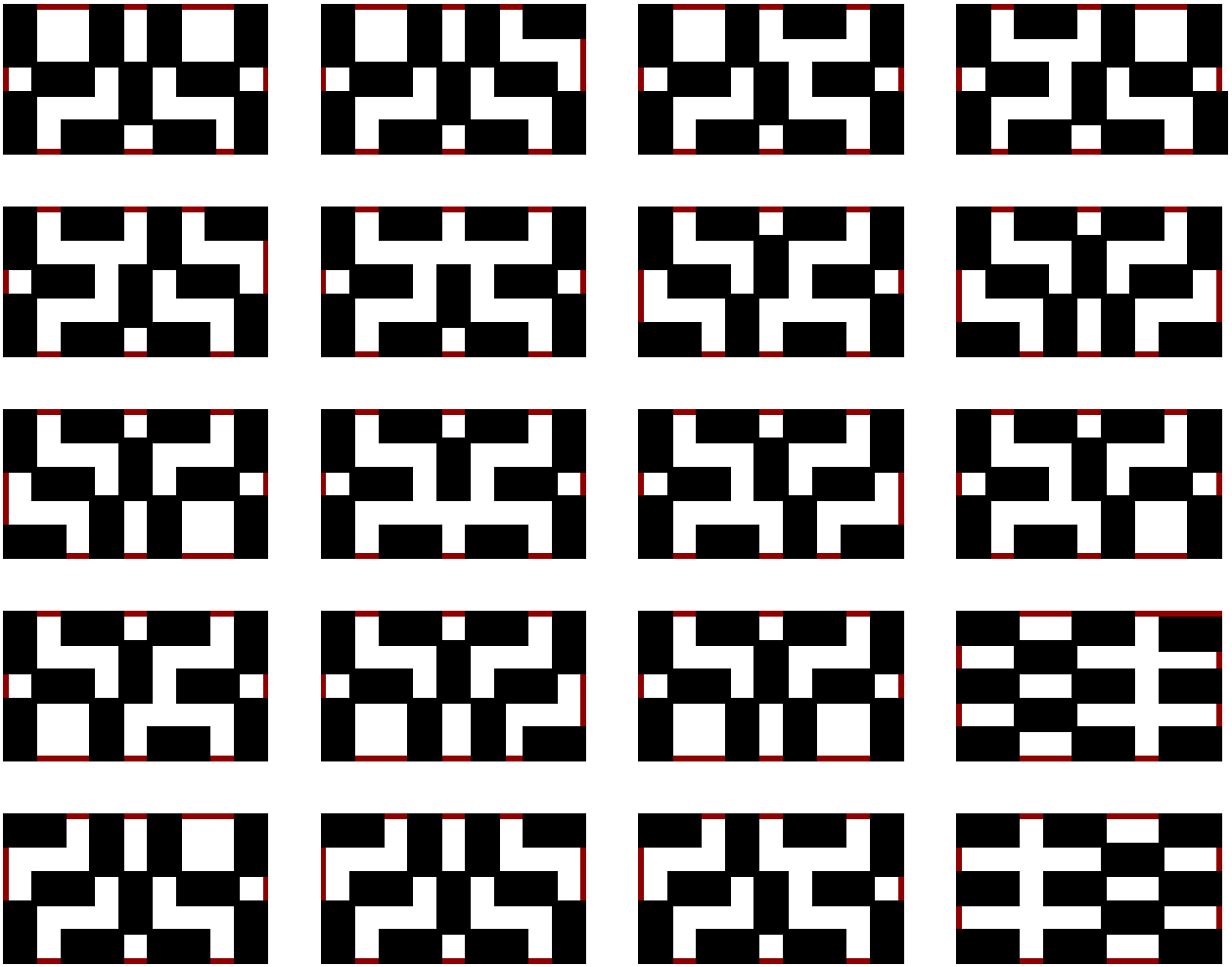}
\caption{The $D(9,5,11)=20$ configurations of 11 non-bonding dominoes on the $9\times 5$ board.
}
\label{Fig.9511}
\end{figure}

If no domino is placed, the empty board is the only solution:
\begin{equation}
D(r,c,0)=1.
\label{eq.Dd0}
\end{equation}
A single domino can be placed on any square, because the
constraint on neighbors is not effective. 
Since we do not pay attention to 
dihedral group symmetries induced by
rotations or flips of the entire configuration
along board middle axes or diagonals or the center, this can be done 
in $r(c-1)$ ways horizontally
plus in $c(r-1)$ ways vertically:
\begin{equation}
D(r,c,1)=2rc-r-c.
\label{eq.Dd1}
\end{equation}

\begin{defn}(Maximum Density)
$\bar d$ is the maximum number of non-bonding dominoes that
can be put on the $r\times c$ board:
\begin{equation}
\bar d(r,c) \equiv \max\{d : D(r,c,d)>0\}.
\label{eq.dbar}
\end{equation}
\end{defn}

Estimates of the maximum number of non-bonding dominoes packed into $r\times c$
rectangles can be derived from Figure \ref{Fig.332}:
\begin{itemize}
\item
If we shift the 2 vertical dominoes in the upper left alignment two squares
down and one square right, assuming a semi-infinite board, a lattice
with a $2\times 2$ unit cell (in the standard definition
of surface physics), one domino per unit cell, arises. The interstitial
voids of non-covered unit squares have shape $2\times 1$.
The asymptotic number of non-bonding dominoes 
proposes $d\sim rc/4$ as an upper limit, but Fig. \ref{Fig.10513}
demonstrates that larger $d$ are possible. 
\item
Another tiling pattern arises if in one of the top middle figures the 2 dominoes
are shifted 3 squares to the right or 3 squares down and cloned.
The upper right constellation in Figure \ref{Fig.747} gives the idea.
(Tiling pattern means: the board can be periodically translated to the right and down, copied,
and the bigger board is still a board of non-bonding dominoes because the dominoes
at the left and right border and top and bottom border of the proto-tile have been
compatible with the adjacency rules.)
This yields a $3\times 3$ unit cell with 2 dominoes; the interstitial
voids are $1\times 1$ and $2\times 2$ squares.
The approximate density is slightly lower, $d\sim 2rc/9$. 
If $r$ and $c$ are both multiples of 3, more than $2rc/9$ dominoes may fit into the rectangle,
demonstrated in Fig.\ \ref{Fig.937}.
\end{itemize}
Edge effects modify these estimates for the finite boards.
\begin{rem}
The figures of the board configurations have been plotted
with \texttt{xfig} based on the sets of coverages
generated by the C++ \texttt{NonAdjDom} program
in the \texttt{anc} directory.
\end{rem}

\begin{defn}(Bivariate GF)
The bivariate generating function for a fixed number of columns of the board is
\begin{equation}
\hat D_c(x,y) = \sum_{r=0}^\infty \sum_{d=0}^\infty D(r,c,d) x^r y^d.
\label{eq.Dhat}
\end{equation}
\end{defn}

\section{Transfer Matrices}
The state of an incomplete board can be described by a 4-ary word
of $c$ letters which characters for all squares in the last row
which dominoes cover any of the $c$ squares. This could be
the alphabet consisting of (i) \texttt{0} which indicates the
square is not covered by a domino, (ii) \texttt{h} which indicates
that the square is covered by a horizontal domino, (iii) \texttt{u} which indicates
the square is covered by a vertical domino which protrudes to the
next row, (iv) \texttt{d} which indicates the square is covered by a vertical
domino that is shared with the previous row. 
By the non-bonding rule, the \texttt{h} must always appear in pairs
of 2, and two letters of \texttt{h,u,l} cannot be adjacent to each
other, which means they must be separated by at least one \texttt{0}.

\begin{figure}
\includegraphics[scale=0.3]{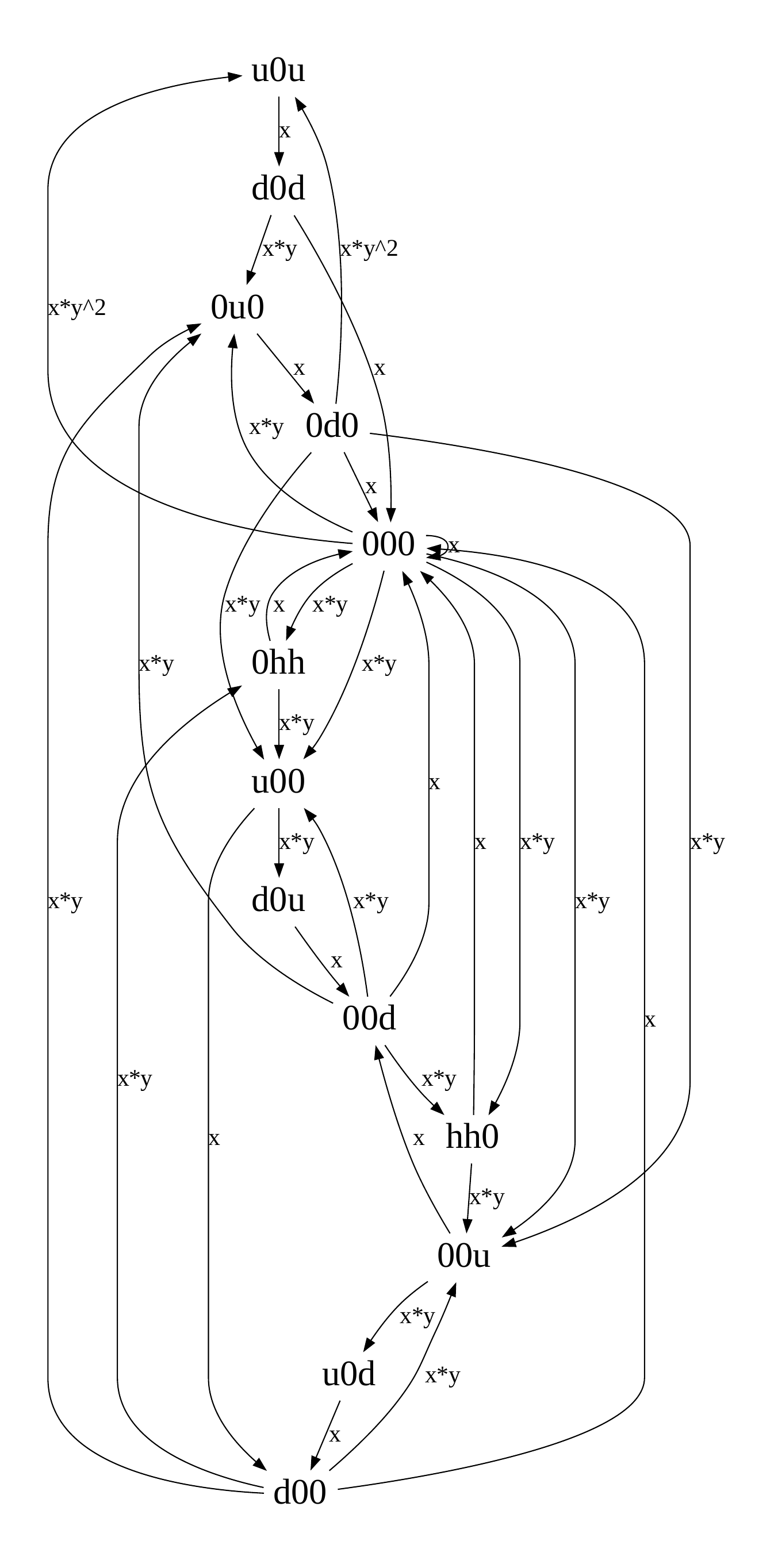}
\caption{State diagram of the $s=13$ states for words with $c=3$ letters.
}
\label{Fig.T3}
\end{figure}

This information is sufficient
to select the set of all words of the next row that are compatible
with the non-bonding rule. These $s$ states are put into a state machine
diagram, like Figure \ref{Fig.T3}. 
The $s\times s$ transfer matrix $T$ is constructed from the transitions
admitted by the no-bonding rule; it contains $xy^w$, where $w$
is the number 
of dominoes added (i.e. the number of \texttt{hh} pairs plus the number of \texttt{u}).
By the Transfer Matrix algorithm the element of the $(1-T)^{-1}$, the
inverse of the unit matrix minus the transfer matrix, which is associated
with the row and column of the all-\texttt{0}-word is the bivariate
generating function.
Compatibility means (i) the next word must be \texttt{d} if the letter
of the previous word is \texttt{u} at the same place, (ii) it must be \texttt{0}
if the letter of the previous word is \texttt{h} or \texttt{d} at
the same place.
\begin{exa}
If $c=1$, the state diagram has $s=3$ vertices: \texttt{d}, \texttt{0}, \texttt{u}.
\end{exa}
\begin{exa}
If $c=2$, the state diagram has $s=6$ vertices: 
\texttt{0d},
\texttt{d0},
\texttt{00},
\texttt{u0},
\texttt{hh},
\texttt{0u}.
\end{exa}
\begin{exa}
If $c=3$, the state diagram, Fig. \ref{Fig.T3}, has $s=13$ vertices: 
\texttt{d0d},
\texttt{00d},
\texttt{u0d},
\texttt{0d0},
\texttt{d00},
\texttt{000},
\texttt{u00},
\texttt{hh0},
\texttt{0u0},
\texttt{0hh},
\texttt{d0u},
\texttt{00u},
\texttt{u0u}.
\end{exa}
\begin{rem}
The number of words grows as $s_{c\ge 0} = 1,3,6,13,28,60,129,277,595,\ldots$
with $s_c=s_{c-1}+2s_{c-2}+s_{c-3}$ \cite[A002478]{sloane}.
The recurrence reflects that a non-bonding word of length $c$ can
be created (i) from any word of length $c-1$ by appending \texttt{0}, (ii) from 
any word of length $c-2$ by appending \texttt{0u} or \texttt{0d}, or (iii) from
any word of length $c-3$ by appending \texttt{0hh}. (This counting argument
resembles a standard counting argument of the Fibonacci sequence.)
\end{rem}
Because only one element of $(1-T)^{-1}$ is wanted, the task is 
to solve a linear system of equations with a $s\times s$ matrix of monomials
in $\{x,y\}$ of that growth  to compute the bivariate GF.

\section{results}

\begin{table}
\begin{tabular}{r|rrrrrrrrrr}
$r\backslash d$ & 0 & 1 & 2 & 3 & 4&5&6 \\
\hline
0&1\\
1&1\\
2&1&1\\
3&1&2\\
4&1&3\\
5&1&4&1\\
6&1&5&3\\
7&1&6&6\\
8&1&7&10&1\\
9&1&8&15&4\\
10&1&9&21&10\\
11&1&10&28&20&1\\
12&1&11&36&35&5&
\end{tabular}
\caption{The number $D(r,1,d)$
of placing $d$ non-bonding dominoes on $r\times 1$ boards \cite[A102547]{sloane}.
Row sums in \cite[A068921]{sloane}.
}
\label{tab.1}
\end{table}

The generating function (GF) associated with Table \ref{tab.1} is
\begin{equation}
\hat D_1(x,y) \equiv p_1(x,y)/q_1(x,y);\quad
p_1(x,y) = 1+x^2y ;\quad q_1(x,y) = 1-x-x^3y.
\label{eq.hatD1}
\end{equation}
The maximum filling is governed by the effective length of the domino,
2 plus the enforced free square, so
\begin{equation}
\bar d(r,1)= \lfloor \frac{r+1}{3}\rfloor.
\end{equation}

\begin{table}
\begin{tabular}{r|rrrrrrrrrrrrrr}
$r\backslash d$ & 0 & 1 & 2 & 3 & 4 &5&6\\
\hline
0&1\\
1&1&1\\
2&1&4\\
3&1&7&1\\
4&1&10&9\\
5&1&13&26&1\\
6&1&16&52&16\\
7&1&19&87&70&1\\
8&1&22&131&190&25\\
9&1&25&184&403&155&1\\
10&1&28&246&736&553&36\\
11&1&31&317&1216&1462&301&1\\
12&1&34&397&1870&3206&1372&49&
\end{tabular}
\caption{The number $D(r,2,d)$
of placing $d$ non-bonding dominoes on $r\times 2$ boards.
Column $d=2$ is \cite[A081267]{sloane}.
}
\label{tab.2}
\end{table}

The GF associated with Table \ref{tab.2} is
\begin{multline}
\hat D_2(x,y) \equiv p_2(x,y)/q_2(x,y);\quad
p_2(x,y) = 1+xy+x^2y-x^3y^2;\\
q_2(x,y) = 1-x-2x^2y-x^3y+x^4y^2.
\label{eq.hatD2}
\end{multline}
The maximum filling is obtained by orienting the long
edge of all dominoes with the short board edge:
\begin{equation}
\bar d(r,2) = \lfloor \frac{r+1}{2} \rfloor.
\end{equation}
The Taylor expansion with respect to $y$ starts
\begin{multline}
\hat D_2(x,y) = \frac{1}{1-x}+x\frac{1+2x}{(1-x)^2}y
+x^3\frac{1+6x+2x^2}{(1-x)^3}y^2
+x^5\frac{1+12x+12x^2+2x^3}{(1-x)^4}y^3+\dots
\end{multline}
\begin{rem}
Coefficients 
$[y^0]\hat D_c(x,y)$
and
$[y^1]\hat D_c(x,y)$
merely rephrase \eqref{eq.Dd0} and \eqref{eq.Dd1}.
\end{rem}
The denominators of the univariate coefficients for $[y^2]\hat D_2$
and $[y^3]\hat D_2$ have the format
$(1-x)^k$, which shows that $D(r,2,2)$ and $D(r,2,3)$ columns in Table \ref{tab.2} 
are polynomials of $r$:
\begin{thm} \cite{Wilf,BernsteinLAA226}\cite[(1.3)]{Gould1972}
\begin{equation}
f(x)\equiv \sum_{r\ge 0}a_rx^r = \frac{\sum_l \gamma_lx^l }{(1-x)^k}
\leftrightarrow
a_r = \sum_l \gamma_l \binom{r+k-l-1}{k-1}.
\end{equation}
\end{thm}
The GF for the row sums is \cite[A105309]{sloane}
\begin{equation}
\hat D_2(x,1) = \frac{1+x+x^2-x^3}{1-x-2x^2-x^3+x^4}.
\end{equation}

\begin{table}
\begin{tabular}{r|rrrrrrrrrrrrr}
$r\backslash d$ & 0 & 1 & 2 & 3 & 4 &5&6&7&8\\
\hline
0&1\\
1&1&2\\
2&1&7&1\\
3&1&12&12\\
4&1&17&45&12\\
5&1&22&103&84&3\\
6&1&27&186&314&92&1\\
7&1&32&294&824&590&60\\
8&1&37&427&1739&2264&726&25\\
9&1&42&585&3184&6467&4234&616&4\\
10&1&47&768&5284&15174&16587&5650&355&1\\
11&1&52&976&8164&30985&50342&30982&5544&149\\
12&1&57&1209&11949&57125&127684&123006&43638&4051&39
\end{tabular}
\caption{The number $D(r,3,d)$
of placing $d$ non-bonding dominoes on $r\times 3$ boards.
$D(3,3,2)=12$ is illustrated in Fig.\ \ref{Fig.332}.
$D(4,3,3)=12$ is illustrated in Fig.\ \ref{Fig.433}.
$D(5,3,4)=3$ is illustrated in Fig.\ \ref{Fig.534}.
}
\label{tab.3}
\end{table}

The GF associated with Table \ref{tab.3} is
\begin{multline}
\hat D_3(x,y) \equiv p_3(x,y)/q_3(x,y);\\
p_3(x,y) = -{x}^{6}{y}^{5}-4{x}^{5}{y}^{4}+{x}^{4}{y}^{3}+2{x}^{3}{y}^{2}+{x}^
{2}{y}^{2}+2xy+{x}^{8}{y}^{6}-2{x}^{6}{y}^{4}+{x}^{4}{y}^{2}+2{x}
^{2}y+1;\\
q_3(x,y) = 
-4{x}^{4}{y}^{2}-2{x}^{4}{y}^{3}-{x}^{9}{y}^{6}+5{x}^{6}{y}^{4}+{
x}^{8}{y}^{6}-3{x}^{5}{y}^{3}+3{x}^{7}{y}^{5}+2{x}^{7}{y}^{4}
\\-{x}
^{5}{y}^{2}-3{x}^{3}{y}^{2}-2{x}^{3}y-x+1-3{x}^{2}y
.
\label{eq.hatD3}
\end{multline}
The rational polynomials of $x$ and $y$ in the
GF's are lengthy for larger $c$. Concise notation
tabulates the coefficients $\alpha$ and $\beta$ in
numerator and denominator:
\begin{defn}(Coefficients of rational GF)
\begin{equation}
\hat D_c(x,y) \equiv \frac{\sum_{i,j \ge 0} \alpha_{c,i,j}x^iy^j}{\sum_{i,j\ge 0} \beta_{c,i,j}x^iy^j}
.
\end{equation}
with the convention that signs are toggled to yield $\alpha_{c,0,0}=\beta_{c,0,0}= +1$.
\end{defn}
With this convention \eqref{eq.hatD3} could also be rewritten as in Table \ref{tab.D3}.
\begin{table}
\begin{tabular}{r|rrrrrrrrrr}
$i\backslash j$ & 0 & 1 & 2 & 3 & 4 & 5 &6\\
\hline
0 & 1 & \\
1 & 0 & 2 & \\
2 & 0 & 2 & 1 & \\
3 & 0 & 0 & 2 & \\
4 & 0 & 0 & 1 & 1 & \\
5 & 0 & 0 & 0 & 0 & -4 & \\
6 & 0 & 0 & 0 & 0 & -2 & -1 & \\
7 & \\
8 & 0 & 0 & 0 & 0 & 0 & 0 & 1 & \\
\end{tabular}
\begin{tabular}{r|rrrrrrrrrr}
$i\backslash j$ & 0 & 1 & 2 & 3 & 4 & 5 &6\\
\hline
0 & 1  &\\
1 & -1  &\\
2 & 0  &-3  &\\
3 & 0  &-2  &-3  &\\
4 & 0  &0  &-4  &-2  &\\
5 & 0  &0  &-1  &-3  &\\
6 & 0  &0  &0  &0  &5  &\\
7 & 0  &0  &0  &0  &2  &3  &\\
8 & 0  &0  &0  &0  &0  &0  &1  &\\
9 & 0  &0  &0  &0  &0  &0  &-1  &\\
\end{tabular} 
\caption{The polynomial
coefficients $\alpha_{3,i,j}$ (left) and $\beta_{3,i,j}$ (right) for the rational GF $\hat D_3(x,y)$.
All $\alpha_{3,7,j}=0$ and this line is empy.
}
\label{tab.D3}
\end{table} 
The nonzero coefficients $\alpha_{c,i,j}$ and $\beta_{c,i,j}$ have
also been collected in files named \texttt{gf}$c$ in the \texttt{anc} directory.
Each line of the file contains 5 fields, blank separated:
\begin{itemize}
\item \texttt{a} if the last integer is $\alpha_{c,i,j}$, \texttt{b} if this is $\beta_{c,i,j}$
\item $c$, the number of columns in the board
\item $i$, the exponent of $x^i$ in the term $x^iy^j$
\item $j$, the exponent of $y^j$ in the term $x^iy^j$
\item $\alpha_{c,i,j}$ or $\beta_{c,i,j}$ depending on the initial letter.
\end{itemize}

The Taylor expansion with respect to $y$ starts
\begin{multline}
\hat D_3(x,y) = \frac{1}{1-x}+x\frac{2+3x}{(1-x)^2}y
+x^2\frac{1+9x+12x^2+3x^3}{(1-x)^3}y^2
\\
+x^4\frac{12+36x+50x^2+24x^3+3x^4}{(1-x)^4}y^3+\dots
\end{multline}
The GF for the row sums is
\begin{equation}
\hat D_3(x,1) = \frac{1+2x+3x^2+2x^3+2x^4-4x^5-3x^6+x^8}{1-x-3x^2-5x^3-6x^4-4x^5+5x^6+5x^7+x^8-x^9},
\end{equation}
obviously obtained by computing row sums of Table \ref{tab.D3}.

\begin{table}
\begin{tabular}{r|rrrrrrrrrrrrr}
$r\backslash d$ & 0 & 1 & 2 & 3 & 4 &5&6&7&8&9&10\\
\hline
0&1\\
1&1&3\\
2&1&10&9\\
3&1&17&45&12\\
4&1&24&126&148&15\\
5&1&31&256&629&349&17\\
6&1&38&435&1758&2327&730&22\\
7&1&45&663&3874&8945&7026&1240&9\\
8&1&52&940&7320&25312&36304&17782&1904&25\\
9&1&59&1266&12439&58880&130822&123240&39512&2799&14\\
10&1&66&1641&19574&119498&372564&561349&361220&78445&3586&17\\
\end{tabular}
\caption{The number $D(r,4,d)$
of placing $d$ non-bonding dominoes on $r\times 4$ boards.
}
\label{tab.4}
\end{table}

The GF $\hat D_4(x,y)$ associated with Table \ref{tab.4} 
has numerator and denominator coefficients of Table \ref{tab.D4}.
\begin{table}
\begin{tabular}{r|rrrrrrrrrrrrrrr}
$i\backslash j$ & 0 & 1 & 2 & 3 & 4 & 5 &6&7&8&9&10&11&12\\
\hline
0 & 1 & \\
1 & 0 & 4 & \\
2 & 0 & 2 & 12 & \\
3 & 0 & 0 & 5 & 20 & \\
4 & 0 & 0 & 1 & 9 & 23 & \\
5 & 0 & 0 & 0 & 1 & -7 & 20 & \\
6 & 0 & 0 & 0 & 0 & -3 & -11 & 18 & \\
7 & 0 & 0 & 0 & 0 & 0 & -3 & 6 & 4 & \\
8 & 0 & 0 & 0 & 0 & 0 & 0 & 3 & 1 & 2 & \\
9 & 0 & 0 & 0 & 0 & 0 & 0 & 0 & 3 & 3 & \\
10 & 0 & 0 & 0 & 0 & 0 & 0 & 0 & 0 & -1 & \\
11 & 0 & 0 & 0 & 0 & 0 & 0 & 0 & 0 & 0 & -1 & -6 & \\
12 & 0 & 0 & 0 & 0 & 0 & 0 & 0 & 0 & 0 & 0 & 0 & -1 & \\
13 & 0 & 0 & 0 & 0 & 0 & 0 & 0 & 0 & 0 & 0 & 0 & 0 & -1 & \\
\end{tabular}
\begin{tabular}{r|rrrrrrrrrrrrrrr}
$i\backslash j$ & 0 & 1 & 2 & 3 & 4 & 5 &6&7&8&9&10&11&12\\
\hline
0 & 1  &\\
1 & -1  &1  &\\
2 & 0  &-6  &\\
3 & 0  &-2  &-23  &-1  &\\
4 & 0  &0  &-8  &-48  &-1  &\\
5 & 0  &0  &-1  &-15  &-55  &\\
6 & 0  &0  &0  &-1  &7  &-55  &\\
7 & 0  &0  &0  &0  &3  &20  &-66  &\\
8 & 0  &0  &0  &0  &0  &3  &3  &-41  &\\
9 & 0  &0  &0  &0  &0  &0  &-3  &-1  &-18  &\\
10 & 0  &0  &0  &0  &0  &0  &0  &-3  &-9  &-13  &\\
11 & 0  &0  &0  &0  &0  &0  &0  &0  &1  &-3  &-9  &\\
12 & 0  &0  &0  &0  &0  &0  &0  &0  &0  &1  &6  &-2  &\\
13 & 0  &0  &0  &0  &0  &0  &0  &0  &0  &0  &0  &1  &-1  &\\
14 & 0  &0  &0  &0  &0  &0  &0  &0  &0  &0  &0  &0  &1  &\\
\end{tabular} 
\caption{The polynomial
coefficients $\alpha_{4,i,j}$ (top) and $\beta_{4,i,j}$ (bottom) for the rational GF $\hat D_4(x,y)$.}
\label{tab.D4}
\end{table} 
The Taylor expansion with respect to $y$ starts
\begin{multline}
\hat D_4(x,y) = \frac{1}{1-x}+x\frac{3+4x}{(1-x)^2}y
+x^2\frac{9+18x+18x^2+4x^3}{(1-x)^3}y^2
\\
+x^3\frac{12+100x+109x^2+82x^3+36x^4+4x^5}{(1-x)^4}y^3+\dots
\end{multline}
The GF for the row sums of Table \ref{tab.4} is
obtained by gathering row sums in Table of \ref{tab.D4},
row sums of $\sum_j \alpha_{4,i,j}$ for the numerator and
$\sum_j \beta_{4,i,j}$ for the denominator of the GF:
\begin{multline}
\hat D_4(x,1) = (
1+4x+14x^2+25x^3+33x^4+14x^5+4x^6+7x^7+6x^8+6x^9-x^{10}
-7x^{11}\\
-x^{12}-x^{13}
)/(
1-x^2-26x^3-57x^4-71x^5-49x^6-43x^6-35x^8-22x^9-25x^{10}\\
-11x^{11}+5x^{12}+x^{14}
)
.
\end{multline}

\begin{table}
\begin{tabular}{r|rrrrrrrrrrrrr}
$r\backslash d$ & 0 & 1 & 2 & 3 & 4 &5&6&7&8\\
\hline
0&1\\
1&1&4&1\\
2&1&13&26&1\\
3&1&22&103&84&3\\
4&1&31&256&629&349&17\\
5&1&40&490&2204&3337&1244&42\\
6&1&49&805&5485&15504&16072&4555&132&1\\
7&1&58&1201&11196&48977&95706&72644&16110&570\\
8&1&67&1678&20066&122063&373374&535434&313802&58927\\
9&1&76&2236&32824&259553&1116890&2528682&2790146&1313712\\
\end{tabular}
\begin{tabular}{r|rrrrrrrrrrrrr}
$r\backslash d$ & 9&10&11\\
\hline
8&1&1854&12\\
9&1&206954&6300&20
\end{tabular}
\caption{The number $D(r,5,d)$
of placing $d$ non-bonding dominoes on $r\times 5$ boards.
}
\label{tab.5}
\end{table}

The GF $\hat D_5(x,y)$ associated with Table \ref{tab.5} 
has the numerator and denominator coefficients 
gathered in the file \texttt{gf5} in the \texttt{anc} directory.

The Taylor expansion with respect to $y$ starts
\begin{multline}
\hat D_5(x,y) = \frac{1}{1-x}+x\frac{4+5x}{(1-x)^2}y
+x\frac{1+23x+28x^2+24x^3+5x^4}{(1-x)^3}y^2
\\
+x^2\frac{1+80x+299x^2+188x^3+108x^4+48x^5+5x^6}{(1-x)^4}y^3+\dots
\end{multline}

The GF $\hat D_6(x,y)$ associated with Table \ref{tab.6} 
has the numerator and denominator coefficients 
gathered in the file \texttt{gf6} in the \texttt{anc} directory.
The Taylor expansion with respect to $y$ starts
\begin{multline}
\hat D_6(x,y) = \frac{1}{1-x}+x\frac{4+5x}{(1-x)^2}y
+x\frac{1+23x+28x^2+24x^3+5x^4}{(1-x)^3}y^2
\\
+x^2\frac{1+80x+299x^2+188x^3+108x^4+48x^5+5x^6}{(1-x)^4}y^3+\dots
\end{multline}

\begin{table}
\begin{tabular}{r|rrrrrrrrrrrrr}
$r\backslash d$ & 0 & 1 & 2 & 3 & 4 &5&6&7&8\\
\hline
0&1&\\
1&1&5&3\\
2&1&16&52&16\\
3&1&27&186&314&92&1\\
4&1&38&435&1758&2327&730&22\\
5&1&49&805&5485&15504&16072&4555&132&1\\
6&1&60&1296&12760&59806&128236&112704&32308&2324\\
7&1&71&1908&24908&168949&593937&1029205&792701&227482\\
8&1&82&2641&43260&390830&1982122&5530922&8055060&5586358\\
\end{tabular}
\begin{tabular}{r|rrrrrrrrrrrrr}
$r\backslash d$ & 9&10&11&12\\
\hline
6&16\\
7&17389&159&1\\
8&1577984&141496&2686&18\\
\end{tabular}
\caption{The number $D(r,6,d)$
of placing $d$ non-bonding dominoes on $r\times 6$ boards.
}
\label{tab.6}
\end{table}

\begin{table}
\begin{tabular}{r|rrrrrrrrrrrrr}
$r\backslash d$ & 0 & 1 & 2 & 3 & 4 &5&6&7&8\\
\hline
0&1\\
1&1&6&6\\
2&1&19&87&70&1\\
3&1&32&294&824&590&60\\
4&1&45&663&3874&8945&7026&1240&9\\
5&1&58&1201&11196&48977&95706&72644&16110&570\\
6&1&71&1908&24908&168949&593937&1029205&792701&227482\\
7&1&84&2784&47200&444544&2373216&7067000&11186540&8632116\\
8&1&97&3829&80269&979481&7189142&31790635&82896379&121946290
\end{tabular}

\begin{tabular}{r|rrrrrrrrrrrrr}
$r\backslash d$ & 9&10&11&12&13&14\\
\hline
6&17389&159&1\\
7&2805332&291836&4708&7\\
8&94134303&34072987&4847487&211850&1550&8
\end{tabular}
\caption{The number $D(r,7,d)$
of placing $d$ non-bonding dominoes on $r\times 7$ boards. The case $D(7,4,7)=D(4,7,7)=9$ is Fig. \ref{Fig.747}.
}
\label{tab.7}
\end{table}

\begin{table}
\begin{tabular}{r|rrrrrrrrrrrrr}
$r\backslash d$ & 0 & 1 & 2 & 3 & 4 &5&6&7\\
\hline
0&1\\
1&1&7&10&1\\
2&1&22&131&190&25\\
3&1&37&427&1739&2264&726&25\\
4&1&52&940&7320&25312&36304&17782&1904\\
5&1&67&1678&20066&122063&373374&535434&313802\\
6&1&82&2641&43260&390830&1982122&5530922&8055060\\
7&1&97&3829&80269&979481&7189142 &31790635 &82896379 \\
8 & 1 & 112 & 5242 & 134468 & 2085933 & 20397928 & 127424616 &  505799712 \\
9  &1 &127 &6880 &209232 &3958594 &48840665 &401073664 &2201053787 \\
\end{tabular}

\begin{tabular}{r|rrrrrrrrrrrrrrrrr}
$r\backslash d$ & 8&9&10&11&12\\
\hline
4&25\\
5&58927 & 1854&12\\
6& 5586358&1577984&141496&2686&18\\
7& 121946290&94134303 & 34072987 & 4847487 & 211850 \\
8& 1249225399&1850015072 & 1553518376 & 682237872 & 139548569 \\
9& 8005746668& 18927854296 &28198554465 &25310043984 &12852685831 \\
\end{tabular}

\begin{tabular}{r|rrrrrrrrrrrrrrrrr}
$r\backslash d$ & 13 & 14 & 15&16&17&18\\
\hline
7 & 1550&8\\
8 & 11218612 &294974& 2112 &16\\
9&  3379935390 &405893730& 18447965 &244776 &1049 &9\\
\end{tabular}
\caption{The number $D(r,8,d)$
of placing $d$ non-bonding dominoes on $r\times 8$ boards.
Row $r=8$ is the relevant entry for chess boards.
}
\label{tab.8}
\end{table}

\begin{table}
\begin{tabular}{r|rrrrrrrrrrrrr}
$r\backslash d$ & 0 & 1 & 2 & 3 & 4 &5&6&7\\
\hline
0&1\\
1&1&8&15&4\\
2&1& 25 & 184 & 403 & 155 & 1\\
3&1& 42 & 585 & 3184 & 6467 & 4234 & 616 & 4\\
4&1& 59 & 1266 & 12439 & 58880 & 130822 & 123240 & 39512 \\
5&1& 76 & 2236 & 32824 & 259553 & 1116890 & 2528682 & 2790146 \\
6&1& 93 & 3495 & 69147 & 787891&  5327437 & 21321326 & 49165252 \\
7&1& 110 & 5043 & 126312 & 1905925 & 18044808 & 108514562& 411721872 \\
8 & 1 &127 &6880 &209232 &3958594 &48840665 &401073664 &2201053787 \\
9 & 1 &144  & 9006 & 322820 & 7374205 & 112989064 & 1191296018 & 8734437124\\
\end{tabular}

\begin{tabular}{r|rrrrrrrrrrrrrrrrr}
$r\backslash d$ & 8&9&10&11&12\\
\hline
4& 2799& 14\\
5& 1313712& 206954 & 6300 & 20\\
6& 62051636& 39434252 & 11086229 & 1145363 & 31238 \\
7& 963380633& 1335889176 & 1033075988& 407606718& 71853631\\
8& 8005746668& 18927854296 &28198554465 &25310043984 &12852685831 \\
9& \ldots\\
\end{tabular}

\begin{tabular}{r|rrrrrrrrrrrrrrrrr}
$r\backslash d$ & 13 & 14 & 15&16&17&18\\
\hline
6 & 100 & 1\\
7 & 4536260& 65621& 110\\
8 & 3379935390 &405893730& 18447965 &244776 &1049 &9\\
9 & \ldots & & & & & 617404\\
\end{tabular}
\caption{The number $D(r,9,d)$
of placing $d$ non-bonding dominoes on $r\times 9$ boards.
Row $r=9$ is incomplete.
}
\label{tab.9}
\end{table}

\begin{rem}
The tables $D(r,c,d)$ presented here are essentially 50\% redundant, implied
by \eqref{eq.Dsym}: 
Row $r=7$ in Table \ref{tab.8} is row $r=8$ in Table \ref{tab.7};
Row $r=6$ in Table \ref{tab.8} is row $r=8$ in Table \ref{tab.6};
Row $r=6$ in Table \ref{tab.7} is row $r=7$ in Table \ref{tab.6};
Row $r=5$ in Table \ref{tab.7} is row $r=7$ in Table \ref{tab.5}; and so on.
\end{rem}

In the data cube of the $D(r,c,d)$ one can also construct other
slices. To focus on boards of square shape, one may gather
row 1 of Table \ref{tab.1},
row 2 of Table \ref{tab.2},
row 3 of Table \ref{tab.3},
row 4 of Table \ref{tab.4}, etc.\ to assemble Table \ref{tab.dia}.

\begin{table}
\begin{tabular}{r|rrrrrrrrrrrrr}
$r\backslash d$ & 0 & 1 & 2 & 3 & 4 &5&6&7&8&9\\
\hline
0&1\\
1&1\\
2&1&4\\
3&1&12&12\\
4&1&24&126&148&15\\
5&1&40&490&2204&3337&1244&42\\
6&1&60&1296&12760&59806&128236&112704&32308&2324&16\\
\end{tabular}
\caption{The number $D(r,r,d)$
of placing $d$ non-bonding dominoes on $r\times r$ boards.
}
\label{tab.dia}
\end{table}

By observing that the $D(r,c,2)$ are apparently all quadratic polynomials in $r$ for $r\ge 3$,
explicitly shown via $[y^2]\hat D_c(x,y)$ for $2\le c\le 6$ above,
and that $D$ must be symmetric observing \eqref{eq.Dsym}, a biquadratic fitting ansatz yields:
\begin{conj}
\begin{equation}
D(r,c,2) = 2c^2r^2-2(cr^2+c^2r)+\frac12(r^2+c^2)-22cr +\frac{59}{2}(c+r)-30;\quad r,c\ge 3.
\end{equation}
\end{conj}
This can be done by taking the quadratic polynomials in $r$, and then constructing
the polynomial in $c$ by Lagrange interpolation.
\begin{rem}
A proof that the bivariate rational GF reduces to polynomials in $r$ could
be founded on the structure of the determinant of $(1-T)^{-1}$ \cite{MatharVixra2404}.
\end{rem}

The statistics of the maximum filling obtained so far
is summarized in Table \ref{tab.dbar}.
\begin{table}
\begin{tabular}{r|rrrrrrrrrr}
$r\backslash c$ & 1 & 2 & 3 & 4 & 5 & 6 & 7 & 8\\
\hline
1& 0 \\
2& 1 & 1\\
3& 1 & 2 & 2\\
4& 1 & 2 & 3 & 4 \\
5& 2 & 3 & 4 & 5 & 6\\
6& 2 & 3 & 5 & 6 & 7 & 9 \\
7& 2 & 4 & 5 & 7 & 8 & 10 & 12 \\
8& 3 & 4 & 6 & 8 & 10 & 12 & 14 & 16 \\
9& 3 & 5 & 7 & 9 & 11 & 14 & 15 & 18\\
10& 3 & 5 & 8 & 10 & 13 & 15 & 18\\
11& 4 & 6 & 8 & 11 & 13 & 17 \\
12& 4 & 6 & 9 & 12 & 15 & 18 \\
\end{tabular}
\caption{Maximum filling $\bar d$, row lengths of the $D(r,c,d)$ tables.
The array is symmetric along
the diagonal, $\bar d(r,c)=\bar d(c,r)$.}
\label{tab.dbar}
\end{table}
Considering Fig.\ \ref{Fig.10513} we find:
\begin{conj}
If either $r$ or $c$ is an odd multiple of 2 and the other odd, $\bar d(r,c)=rc/4+1/2$.
\end{conj}
\begin{conj}
If $r$ and $c$ are both even, $\bar d(r,c)=rc/4$.
\end{conj}
\begin{conj}
If either $r$ or $c$ is a multiple of 4 and the other odd $\ge 3$, $\bar d(r,c)=rc/4$.
\end{conj}

\section{Summary} 
The bivariate GF \eqref{eq.Dhat} which enumerates $r \times c$ boards
with $d$ non-bonding non-overlapping dominoes has been presented for $1\le c\le 6$
in Eqs.\ 
\eqref{eq.hatD1},
\eqref{eq.hatD2},
\eqref{eq.hatD3},
Table \ref{tab.D4},
and in machine-readable form in the \texttt{anc} directory.

\bibliographystyle{amsplain}
\bibliography{all}

\providecommand{\bysame}{\leavevmode\hbox to3em{\hrulefill}\thinspace}
\providecommand{\MR}{\relax\ifhmode\unskip\space\fi MR }
\providecommand{\MRhref}[2]{%
  \href{http://www.ams.org/mathscinet-getitem?mr=#1}{#2}
}
\providecommand{\href}[2]{#2}
\begin{thebibliography}{1}

\bibitem{BernsteinLAA226}
Mira Bernstein and Neil J.~A. Sloane, \emph{Some canonical sequences of
  integers}, Lin.\ Alg.\ Applic. \textbf{226--228} (1995), 57--72, (E:)
  \cite{BrualdiLAA320}. \MR{1344554}

\bibitem{BrualdiLAA320}
Richard~A. Brualdi, \emph{From the editor-in-chief}, Lin.\ Alg.\ Applic.
  \textbf{320} (2000), no.~1--3, 209--216. \MR{1796542}

\bibitem{sloane}
O.~E. I.~S. Foundation~Inc., \emph{The {O}n-{L}ine {E}ncyclopedia {O}f
  {I}nteger {S}equences},  (2024), https://oeis.org/. \MR{3822822}

\bibitem{Gould1972}
Henry~W. Gould, \emph{Combinatorial identities}, 1972. \MR{0354401}

\bibitem{MatharVixra2404}
Richard~J. Mathar, \emph{Bivariate generating functions for non-attacking
  wazirs on rectangular boards}, viXra:2404.0122 (2024).

\bibitem{SiehlerArxiv1409}
Jacob~A. Siehler, \emph{Selections without adjacency on a rectangular grid},
  arXiv:1409.3869 (2014).

\bibitem{Wilf}
Herbert~S. Wilf, \emph{Generatingfunctionology}, Academic Press, 2004.
  \MR{2172781}

\end{thebibliography}

\end{document}